\documentclass[a4paper]{article}
\usepackage[T1]{fontenc}
\usepackage{url}
\usepackage[dvipdfmx]{graphicx}
\usepackage{latexsym}
\usepackage[numbers]{natbib}
\usepackage{rotating}
\usepackage{amsmath,amssymb,amsfonts}
\usepackage{algorithm,algorithmic}
\usepackage{textcomp}
\usepackage{xcolor}

\begin{document}

\title{Assessing the Performance of Mixed-Precision ILU(0)-Preconditioned Multiple-Precision Real and Complex Krylov Subspace Methods}

\author{Tomonori Kouya\\Orcid ID: 0000-0003-0178-5519\\Otemon Gakuin University}

\maketitle

\begin{abstract}
Krylov subspace methods are linear solvers based on matrix-vector multiplications and vector operations. While easily parallelizable, they are sensitive to rounding errors and may experience convergence issues. ILU(0), an incomplete LU factorization with zero fill-in, is a well-known preconditioning technique that enhances convergence for sparse matrices. In this paper, we implement a double-precision and multiple-precision ILU(0) preconditioner, compatible with product-type Krylov subspace methods, and evaluate its performance.

\end{abstract}

\section{Introduction}

With the continuous progress of digital transformation (DX), computing infrastructures capable of efficiently processing large-scale data are increasingly essential. The operation and training of AI systems, primarily based on machine learning and deep learning, are expected to drive DX forward. Consequently, faster and more efficient computers are needed, and the demand for computational power continues to grow. However, further increases in CPU and GPU clock speeds are limited, making adopting architectures optimized for parallel processing, such as those supporting single instruction, multiple data (SIMD) instructions and multi-core configurations essential. Software should likewise be optimized to run effectively on these modern architectures.

While AI applications frequently use floating-point formats with short mantissas (e.g., IEEE754 binary32 or binary16), scientific computing requires higher precision formats such as binary64 or above. Hardware advancements are driven by market demand, which often favors the former. For ill-conditioned problems requiring high precision or mixed-precision computing that combines short and long mantissa floating-point formats is a promising approach to achieve both accuracy and performance.

To address such ill-conditioned problems, we developed an optimized basic linear algebra library, BNCmatmul, supporting multiple-precision formats beyond binary64. Libraries such as QD \cite{qd}, MPFR \cite{mpfr}, and MPC \cite{mpc} offer reliable multi-component and arbitrary-precision arithmetic. MPLAPACK/MPBLAS \cite{mplapack} extends LAPACK/BLAS functionality to support these formats. However, current parallel architecture optimizations are limited to what is provided by C++ compilers. Similar to machine-optimized libraries like Intel MKL \cite{imkl} and OpenBLAS \cite{openblas}, optimized multiple-precision BLAS libraries are required.

Our BNCmatmul library is optimized primarily for x86 environments. As of March 2025, we have completed the implementation of real and complex various multiple-precision floating-point operations, along with random sparse matrix-vector multiplication (SpMV). Details on SpMV are reported in \cite{kouya_iceet2024}, and although further tuning is possible, the current implementation is usable functionally.

This paper reports on applying ILU(0), factorization-free of fill-ins-as a preconditioner to product-type Krylov subspace methods: BiCG, CGS, BiCGSTAB, and GPBiCG. Unlike our AVX2-optimized SpMV, the current ILU(0) implementation is not optimized or parallelized and is, therefore, slow. Nevertheless, the implementation supports mixed-precision ILU(0) using binary64, enabling performance evaluation with test matrices from the SuiteSparse matrix collection \cite{ufsparse}. This evaluation offers insight into the utility of the method and provides a foundation for future optimization.

The remainder of this paper is structured as follows. Section 2 describes the features of BNCmatmul. Section 3 explains the ILU(0) algorithm. Section 4 presents benchmark test results on an EPYC machine and discussions. Section 5 concludes with future work.

\section{Multiple-Precision Optimized Basic Linear Algebra Library: BNCmatmul}

\begin{figure}[htb]
  \centering
  \includegraphics[width=.99\textwidth]{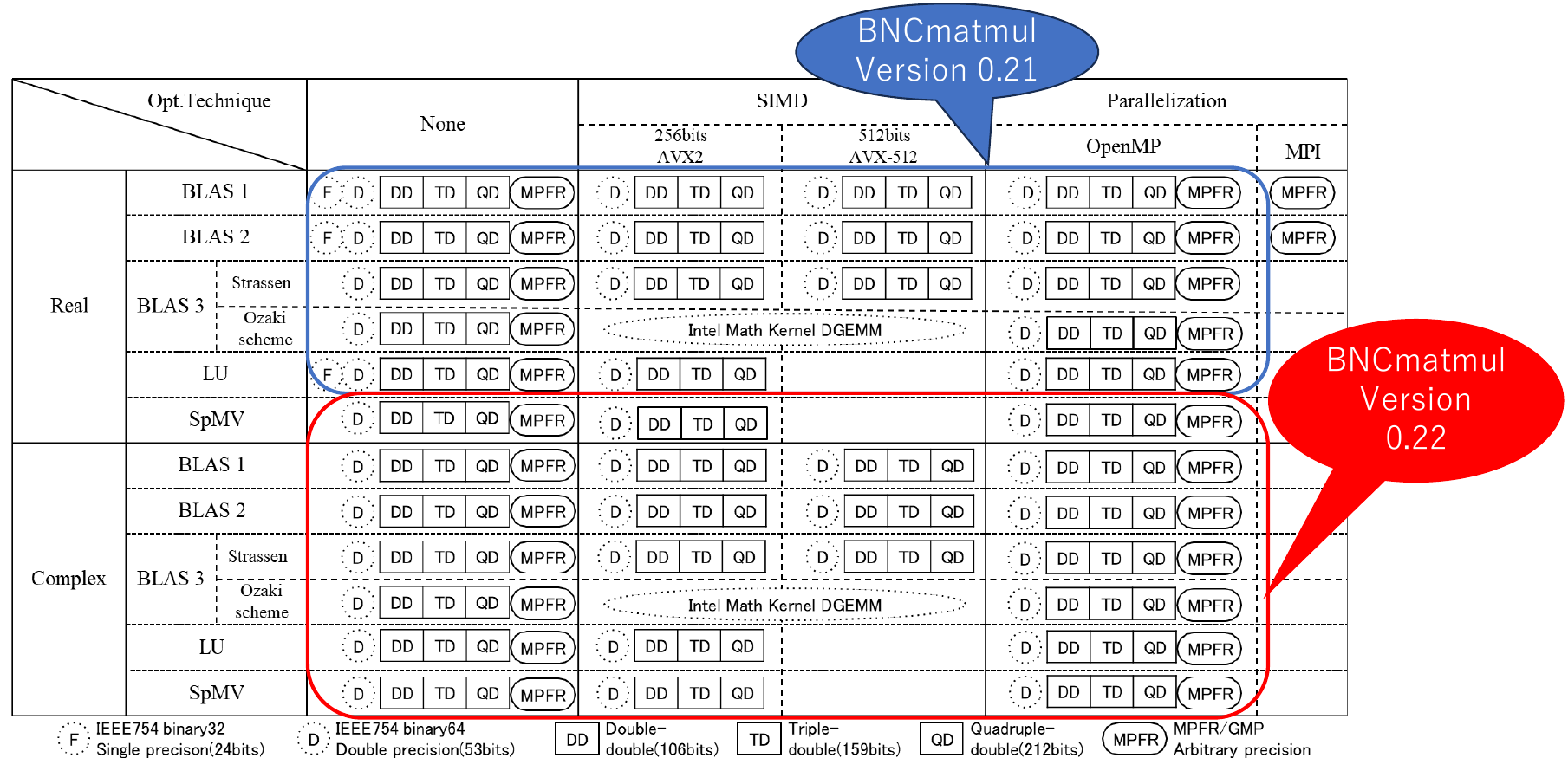}
    \caption{Current development status of BNCmatmul}
  \label{fig:bncmatmul}
\end{figure}

An overview of the features of BNCmatmul \cite{bncmatmul} is shown in Figure~\ref{fig:bncmatmul}. As of March 2025, Version 0.21 has been released as open-source software, primarily supporting real-valued linear algebra routines and their optimized/parallelized versions. Version 0.22, currently under development, extends functionality to complex linear algebra based on real implementations.

The library supports the following five floating-point types:

\begin{description}
  \item[Double (D)] IEEE754 binary64 format with a 53-bit mantissa (52-bit + 1 hidden bit).
  \item[Multi-component-type] Formats based on binary64 that extend mantissa length using multiple components:
  \begin{description}
    \item[Double-double (DD)] 106-bit mantissa using two binary64 values.
    \item[Triple-double (TD)] 159-bit mantissa using three binary64 values.
    \item[Quadruple-double (QD)] 212-bit mantissa using four binary64 values.
  \end{description}
  \item[Multi-digit-type] Arbitrary precision types provided by MPFR (real) and MPC (complex), built on GMP's multiple-precision natural number (MPN) kernel.
\end{description}

The SpMV implemented in this work builds on our prior implementation of MPFR-based arbitrary-precision compressed row storage (CRS) multiplication \cite{kouya_sparse2012}, extended to support all five floating-point types. This feature is included in Version 0.22, and detailed performance evaluation results can be found in \cite{kouya_iceet2024}.

Version 0.22 is scheduled for release soon. Looking ahead, Version 0.23 will expand the capabilities of the library for nonlinear computations based on the current linear algebra core. It will also introduce SIMD-optimized "pair arithmetic" for multi-component types. Additionally, we plan to adapt the library for consumer-grade computing platforms, including environments with 128-bit SIMD support such as Arm NEON, WebAssembly (WASM), and x86 SSE.

\section{Product-Type Krylov Subspace Methods and Mixed-Precision ILU(0) Preconditioner}

We consider solving the linear system:
\begin{equation}
  A\mathbf{x} = \mathbf{b}
  \label{eq:linear_eq}
\end{equation}
where $A \in \mathbb{C}^{n \times n}$ is a complex-valued square coefficient matrix, $\mathbf{b} \in \mathbb{C}^n$ is a known $n$-dimensional vector, and $\mathbf{x} \in \mathbb{C}^n$ is the unknown $n$-dimensional vector to be solved. In this study, we assume $A$ is non-singular. For comparison purposes, we also consider real-valued cases where $A \in \mathbb{R}^{n \times n}$, $\mathbf{b} \in \mathbb{R}^n$, and $\mathbf{x} \in \mathbb{R}^n$.

When $A$ is sparse, iterative methods such as Krylov subspace methods and GMRES are effective due to their reliance on matrix-vector multiplication, which is easy to parallelize. However, these algorithms are highly sensitive to rounding errors and may fail to converge within the theoretical limit of $n$ iterations. To improve convergence, preconditioning is commonly applied. One such technique is the incomplete LU factorization (ILU) that approximates $A \approx LU$ while controlling fill-ins. We use ILU(0), which introduces no fill-in, as defined in the Algorithm~\ref{algo:ilu0}.

\begin{algorithm}
  \caption{ILU(0) Factorization}
  \label{algo:ilu0}
  \begin{algorithmic}[1]
    \REQUIRE Sparse matrix $A \in \mathbb{C}^{n \times n}$ (non-singular)
    \ENSURE Matrices $\tilde{L}$, $\tilde{U} \in \mathbb{C}^{n \times n}$ such that $A \approx \tilde{L}\tilde{U}$
    \STATE Initialize $\tilde{L} := I$, $\tilde{U} := 0$
    \FOR{$k = 1$ to $n$}
     \IF{$U_{kk} := A_{kk} \not= 0$}
      \FOR{each $i > k$ such that $A_{ik} \neq 0$}
        \STATE $\tilde{L}_{ik} := A_{ik} / \tilde{U}_{kk}$
        \FOR{each $j > k$ such that $A_{kj} \neq 0$}
          \IF{$A_{ij} \neq 0$}
            \STATE $A_{ij} := A_{ij} - \tilde{L}_{ik} \cdot \tilde{U}_{kj}$
          \ENDIF
        \ENDFOR
      \ENDFOR
      \FOR{each $j > k$ such that $A_{kj} \neq 0$}
        \STATE $\tilde{U}_{kj} := A_{kj}$
      \ENDFOR
     \ENDIF
    \ENDFOR
  \end{algorithmic}
\end{algorithm}

Once the ILU(0) factorization is complete, we apply the preconditioner $M = \tilde{L}\tilde{U}$ before starting the iterations. During each iteration, we solve $M\mathbf{z} = \mathbf{r}$ using forward and backward substitutions, as described in Algorithm~\ref{algo:ilu0_fbsubst}.

\begin{algorithm}
  \caption{Preconditioning via ILU(0): Solving $M^{-1}\mathbf{r} = \mathbf{z}$}
  \label{algo:ilu0_fbsubst}
  \begin{algorithmic}[1]
    \REQUIRE Residual vector $\mathbf{r} \in \mathbb{C}^n$, ILU(0) factors $\tilde{L}$, $\tilde{U}$
    \ENSURE Preconditioned vector $\mathbf{z} \in \mathbb{C}^n$
    \STATE Forward substitution to solve $\tilde{L}\mathbf{y} = \mathbf{r}$
    \FOR{$i = 1$ to $n$}
      \STATE $y_i := r_i$
      \FOR{$j = 1$ to $i - 1$}
        \IF{$\tilde{L}_{ij} \neq 0$}
          \STATE $y_i := y_i - \tilde{L}_{ij} \cdot y_j$
        \ENDIF
      \ENDFOR
    \ENDFOR
    \STATE Backward substitution to solve $\tilde{U}\mathbf{z} = \mathbf{y}$
    \FOR{$i = n$ downto $1$}
      \STATE $z_i := y_i$
      \FOR{$j = i+1$ to $n$}
        \IF{$\tilde{U}_{ij} \neq 0$}
          \STATE $z_i := z_i - \tilde{U}_{ij} \cdot z_j$
        \ENDIF
      \ENDFOR
      \STATE $z_i := z_i / \tilde{U}_{ii}$
    \ENDFOR
  \end{algorithmic}
\end{algorithm}

Among the product-type Krylov methods, applying ILU(0) preconditioning to the BiCG method results in Algorithm~\ref{algo:precond_bicg}:

\begin{algorithm}
  \caption{Preconditioned BiCG Method for Complex Systems}
  \label{algo:precond_bicg}
  \begin{algorithmic}[1]
    \REQUIRE Complex matrix $A \in \mathbb{C}^{n \times n}$, RHS vector $\mathbf{b} \in \mathbb{C}^n$, initial guess $\mathbf{x}_0 := 0$
    \REQUIRE Preconditioner $M \approx A$, tolerance $\epsilon$
    \ENSURE Approximate solution $\mathbf{x}$
    \STATE $\mathbf{x} := \mathbf{x}_0$
    \STATE $\mathbf{r}_0 := \mathbf{b} - A\mathbf{x}$
    \STATE $\mathbf{r} := \mathbf{r}_0 $
    \STATE $\tilde{\mathbf{r}} := \mathbf{r}$
    \STATE Solve $M \mathbf{z} := \mathbf{r}$
    \STATE Solve $M^H \tilde{\mathbf{z}} := \tilde{\mathbf{r}}$
    \STATE $\mathbf{p} := \mathbf{z}$, $\tilde{\mathbf{p}} := \tilde{\mathbf{z}}$
    \STATE $\rho := \langle \tilde{\mathbf{z}}, \mathbf{r} \rangle$
    \FOR{$k = 1$ to max\_iter}
      \STATE $\mathbf{q} := A\mathbf{p}$
      \STATE $\tilde{\mathbf{q}} := A^H \tilde{\mathbf{p}}$
      \STATE $\alpha := \rho / \langle \tilde{\mathbf{p}}, \mathbf{q} \rangle$
      \STATE $\mathbf{x} := \mathbf{x} + \alpha \mathbf{p}$
      \STATE $\mathbf{r} := \mathbf{r} - \alpha \mathbf{q}$
      \STATE $\tilde{\mathbf{r}} := \tilde{\mathbf{r}} - \overline{\alpha} \tilde{\mathbf{q}}$
      \IF{$\|\mathbf{r}\| < \varepsilon_r \|\mathbf{r}_0\| + \varepsilon_a $}
        \STATE \textbf{break}
      \ENDIF
      \STATE Solve $M\mathbf{z} := \mathbf{r}$
      \STATE Solve $M^H \tilde{\mathbf{z}} := \tilde{\mathbf{r}}$
      \STATE $\rho_{\text{new}} := \langle \tilde{\mathbf{z}}, \mathbf{r} \rangle$
      \STATE $\beta := \rho_{\text{new}} / \rho$
      \STATE $\mathbf{p} := \mathbf{z} + \beta \mathbf{p}$
      \STATE $\tilde{\mathbf{p}} := \tilde{\mathbf{z}} + \overline{\beta} \tilde{\mathbf{p}}$
      \STATE $\rho := \rho_{\text{new}}$
    \ENDFOR
  \end{algorithmic}
\end{algorithm}

In BiCG, both $M\mathbf{z} = \mathbf{r}$ and $M^H \mathbf{z} = \mathbf{r}$ must be solved at each iteration. In contrast, other product-type methods such as CGS, BiCGSTAB, and GPBiCG require only one preconditioning solve per iteration.

Our ILU(0) implementation supports both real and complex CRS sparse matrices and handles Hermitian systems. In line with previous research by Hishinuma et al. \cite{dd_avx_original}, the SuiteSparse Matrix Collection \cite{ufsparse} provides only binary64 matrices, therefore, we implemented mixed-precision SpMV and forward/backward substitutions with binary64 and higher precisions to assess the impact of mixed-precision ILU(0) on convergence and performance.

\section{Benchmark Test}

This section presents and discusses the results of benchmark tests conducted on the following 32-core EPYC computing environment:

\begin{description}
  \setlength{\parskip}{0cm}
  \item[Hardware:] AMD EPYC 9354P, 3.7 GHz, 32 cores, Ubuntu 20.04.6 LTS
  \item[Software:] Intel Compiler version 2021.10.0, GNU MP 6.2.1, MPFR 4.1.0, MPC 1.2.1
\end{description}

For the linear system defined in Equation~\ref{eq:linear_eq}, the exact solution vector $\mathbf{x}$ is set as follows:
\begin{itemize}
  \item For real matrices: $\mathbf{x} := \sqrt{2}\ [1\ 2\ \ldots\ n]^T$
  \item For complex matrices: $\mathbf{x} := \sqrt{2 + 3\mathrm{i}}\ [1\ 2\ \ldots\ n]^T$
\end{itemize}
The right-hand side vector is computed as $\mathbf{b} := A\mathbf{x}$.

\begin{figure}[htb]
  \centering
  \includegraphics[width=.44\textwidth]{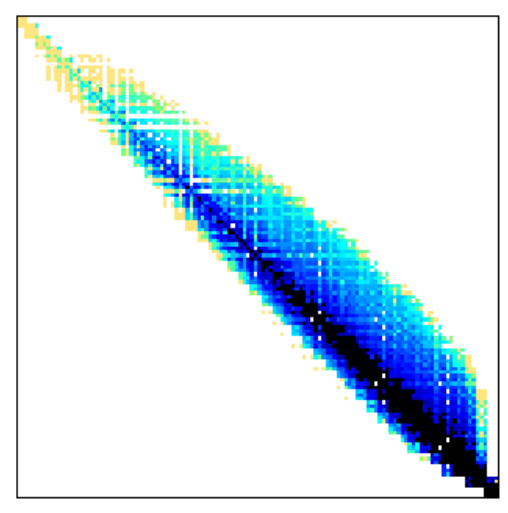}
  \includegraphics[width=.44\textwidth]{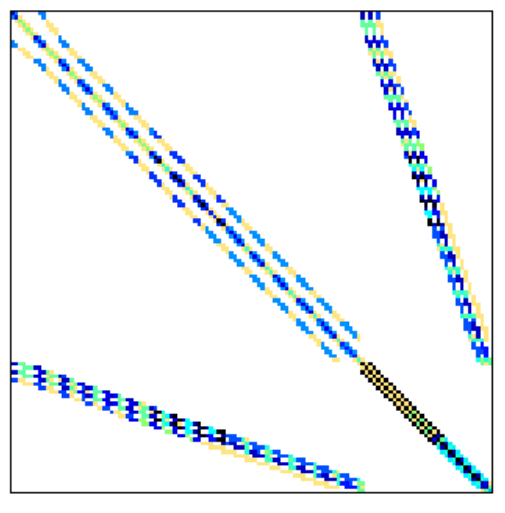}
  \caption{Structures of non-zero elements in the example matrices: \texttt{mpfe} (left) and \texttt{dwg961b} (right) adapted from \cite{ufsparse}}
  \label{fig:example_sparse_matrices}
\end{figure}

The performance of sparse matrix preconditioning strongly depends on the structure of the matrix and the target computing environment. In this test, we selected the following two example matrices from the SuiteSparse Matrix Collection, shown in Figure~\ref{fig:example_sparse_matrices}:
\begin{itemize}
  \item Real sparse matrix: \texttt{mpfe} ($n=768$)
  \item Complex sparse matrix: \texttt{dwg961b} ($n=961$)
\end{itemize}

As demonstrated subsequently, neither matrix converges without ILU(0) preconditioning under floating-point precision below 256 bits. Since both are originally stored in binary64 precision, we evaluated the performance using:
\begin{itemize}
  \item Mixed-precision SpMV implementations using binary64 and DD/TD/QD (BiCG\_d, etc.)
  \item Mixed-precision ILU(0) forward/backward substitution using binary64 (BiCG\_d\_ILU(0), etc.)
  \item Parallel implementations using 32-thread OpenMP without preconditioning (BiCG\_p, etc.)
\end{itemize}

We first present results using DD, TD, and QD multi-component precision formats. Table~\ref{tab:real} lists the results for the real matrix \texttt{mpfe} ($n=765$) under a maximum iteration count of $3n = 2295$, with convergence thresholds $\varepsilon_r = 10^{-13}$ and $\varepsilon_a = 10^{-100}$. The metrics in each column are as follows:
\begin{itemize}
  \item[\textbf{[1]}] Number of iterations
  \item[\textbf{[2]}] Total computation time (in seconds)
  \item[\textbf{[3]}] Average time per iteration (in milliseconds)
\end{itemize}

\begin{table}[htb]
    \caption{Real linear equation: mpfe}
    \label{tab:real}
    \centering
    {\small \begin{tabular}{|c|rrr|rrr|rrr|}\hline
        mcfe, $n=765$ & \multicolumn{3}{c|}{DD} & \multicolumn{3}{c|}{TD} & \multicolumn{3}{c|}{QD} \\ 
        Algorithm & \multicolumn{1}{c|}{[1]} & \multicolumn{1}{c|}{[2]} & \multicolumn{1}{c|}{[2]/[1]} & \multicolumn{1}{c|}{[1]} & \multicolumn{1}{c|}{[2]} & \multicolumn{1}{c|}{[2]/[1]} & \multicolumn{1}{c|}{[1]} & \multicolumn{1}{c|}{[2]} & \multicolumn{1}{c|}{[2]/[1]} \\ \hline
BiCG & 2295 & 0.26 & 0.11 & 2295 & 3.00 & 1.31 & 2295 & 3.64 & 1.58 \\
BiCG\_p & 2295 & 0.23 & 0.10 & 2295 & 1.14 & 0.50 & 2295 & 1.70 & 0.74 \\
BiCG\_ILU(0)  & 13 & 0.67 & 51.31 & 13 & 0.74 & 56.85 & 13 & 0.78 & 59.77 \\
BiCG\_d & 2295 & 0.25 & 0.11 & 2295 & 2.88 & 1.25 & 2295 & 3.47 & 1.51 \\
BiCG\_dp  & 2295 & 0.27 & 0.12 & 2295 & 1.13 & 0.49 & 2295 & 1.61 & 0.70 \\
BiCG\_d\_ILU(0)  & 13 & 0.67 & 51.31 & 13 & 0.69 & 53.23 & 13 & 0.73 & 56.23 \\ \hline
CGS  & 2295 & 0.25 & 0.11 & 2295 & 3.10 & 1.35 & 2295 & 4.01 & 1.75 \\
CGS\_p  & 2295 & 0.09 & 0.04 & 2295 & 0.70 & 0.30 & 2295 & 0.93 & 0.41 \\
CGS\_ILU(0) & 9 & 0.52 & 57.67 & 9 & 0.55 & 60.67 & 9 & 0.59 & 66.00 \\
CGS\_d  & 2295 & 0.23 & 0.10 & 2295 & 2.98 & 1.30 & 2295 & 3.84 & 1.67 \\
CGS\_dp & 2295 & 0.22 & 0.10 & 2295 & 0.69 & 0.30 & 2295 & 0.91 & 0.40 \\
CGS\_d\_ILU(0) & 9 & 0.47 & 52.56 & 9 & 0.51 & 56.56 & 9 & 0.54 & 59.78 \\ \hline
BiCGSTAB & 2295 & 0.25 & 0.11 & 2295 & 3.19 & 1.39 & 2295 & 4.23 & 1.84 \\
BiCGSTAB\_p & 2295 & 0.08 & 0.04 & 2295 & 0.81 & 0.35 & 2295 & 1.12 & 0.49 \\
BiCGSTAB\_ILU(0)  & 9 & 0.51 & 56.89 & 9 & 0.52 & 58.11 & 9 & 0.56 & 62.56 \\
BiCGSTAB\_d & 2295 & 0.24 & 0.10 & 2295 & 3.08 & 1.34 & 2295 & 4.06 & 1.77 \\
BiCGSTAB\_dp  & 2295 & 0.24 & 0.10 & 2295 & 0.83 & 0.36 & 2295 & 1.11 & 0.48 \\
BiCGSTAB\_d\_ILU(0)  & 9 & 0.47 & 52.44 & 9 & 0.49 & 54.11 & 9 & 0.51 & 57.00 \\ \hline
GPBiCG & 2295 & 0.26 & 0.12 & 2295 & 3.53 & 1.54 & 2295 & 4.57 & 1.99 \\
GPBiCG\_p & 2295 & 0.12 & 0.05 & 2295 & 1.20 & 0.52 & 2295 & 1.49 & 0.65 \\
GPBiCG\_ILU(0)  & 8 & 0.46 & 57.63 & 8 & 0.50 & 62.38 & 8 & 0.54 & 67.38 \\
GPBiCG\_d & 2295 & 0.26 & 0.11 & 2295 & 3.42 & 1.49 & 2295 & 4.40 & 1.92 \\
GPBiCG\_dp  & 2295 & 0.23 & 0.10 & 2295 & 1.13 & 0.49 & 2295 & 1.47 & 0.64 \\ 
GPBiCG\_d\_ILU(0)  &8	&0.43	&53.50	&9	&0.51	&56.56	&8	&0.49	&61.13\\ \hline
    \end{tabular}}
\end{table}

As shown, without ILU(0) preconditioning, all methods hit the iteration limit without converging, regardless of the precision format. Mixed-precision SpMV does not negatively impact convergence, nor does it reduce the effectiveness of ILU(0). Mixed-precision forward/backward substitution maintains the same number of iterations.

Next, Table~\ref{tab:complex} lists the results for the complex sparse matrix \texttt{dwg961b} ($n=961$), again with DD, TD, and QD precision formats. As in the previous case, none of the Krylov subspace methods converge without ILU(0) preconditioning.

\begin{table}[htb]
    \caption{Complex linear equation: dwg961b}
    \label{tab:complex}
    \begin{center}
    {\small \begin{tabular}{|c|rrr|rrr|rrr|}\hline
        dwg961b, $n=961$ & \multicolumn{3}{c|}{DD} & \multicolumn{3}{c|}{TD} & \multicolumn{3}{c|}{QD} \\ 
        Algorithm & \multicolumn{1}{c|}{[1]} & \multicolumn{1}{c|}{[2]} & \multicolumn{1}{c|}{[2]/[1]} & \multicolumn{1}{c|}{[1]} & \multicolumn{1}{c|}{[2]} & \multicolumn{1}{c|}{[2]/[1]} & \multicolumn{1}{c|}{[1]} & \multicolumn{1}{c|}{[2]} & \multicolumn{1}{c|}{[2]/[1]} \\ \hline
        BiCG & 2883 & 0.79 & 0.27 & 2883 & 8.96 & 3.11 & 2883 & 13.10 & 4.55 \\
        BiCG\_p & 2883 & 1.61 & 0.56 & 2883 & 6.73 & 2.33 & 2883 & 10.77 & 3.73 \\
        BiCG\_ILU(0)  & 758 & 199.67 & 263.41 & 516 & 170.33 & 330.09 & 442 & 206.52 & 467.24 \\
        BiCG\_d & 2883 & 0.76 & 0.26 & 2883 & 8.64 & 3.00 & 2883 & 12.66 & 4.39 \\
        BiCG\_dp  & 2883 & 1.83 & 0.63 & 2883 & 6.71 & 2.33 & 2883 & 10.36 & 3.59 \\
        BiCG\_d\_ILU(0)  & 758 & 62.74 & 82.77 & 512 & 45.34 & 88.55 & 438 & 82.33 & 187.96 \\ \hline
        CGS  & 2883 & 0.85 & 0.29 & 2883 & 9.89 & 3.43 & 2883 & 15.49 & 5.37 \\
        CGS\_p  & 2883 & 0.70 & 0.24 & 2883 & 3.67 & 1.27 & 2883 & 6.04 & 2.10 \\
        CGS\_ILU(0) & 1057 & 268.15 & 253.69 & 666 & 215.51 & 323.59 & 529 & 244.34 & 461.88 \\
        CGS\_d  & 2883 & 0.81 & 0.28 & 2883 & 9.60 & 3.33 & 2883 & 14.99 & 5.20 \\
        CGS\_dp & 2883 & 1.67 & 0.58 & 2883 & 3.72 & 1.29 & 2883 & 6.07 & 2.10 \\
        CGS\_d\_ILU(0) & 1064 & 84.79 & 79.69 & 664 & 57.05 & 85.92 & 526 & 48.36 & 91.93 \\ \hline
        BiCGSTAB & 2883 & 0.93 & 0.32 & 2883 & 10.51 & 3.65 & 2883 & 17.56 & 6.09 \\
        BiCGSTAB\_p & 2883 & 0.69 & 0.24 & 2883 & 4.27 & 1.48 & 2883 & 8.19 & 2.84 \\
        BiCGSTAB\_ILU(0)  & 2209 & 551.94 & 249.86 & 1242 & 401.85 & 323.55 & 778 & 354.55 & 455.72 \\
        BiCGSTAB\_d & 2883 & 0.89 & 0.31 & 2883 & 10.21 & 3.54 & 2883 & 17.10 & 5.93 \\
        BiCGSTAB\_dp  & 2883 & 1.50 & 0.52 & 2883 & 4.23 & 1.47 & 2883 & 8.17 & 2.83 \\
        BiCGSTAB\_d\_ILU(0)  & 1998 & 159.13 & 79.65 & 1117 & 96.01 & 85.96 & 831 & 76.92 & 92.56 \\ \hline
        GPBiCG & 2883 & 1.12 & 0.39 & 2883 & 12.53 & 4.35 & 2883 & 21.92 & 7.60 \\
        GPBiCG\_p & 2883 & 0.92 & 0.32 & 2883 & 6.30 & 2.18 & 2883 & 12.64 & 4.38 \\
        GPBiCG\_ILU(0)  & 1954 & 489.50 & 250.51 & 1143 & 370.68 & 324.31 & 869 & 396.26 & 456.00 \\
        GPBiCG\_d & 2883 & 1.08 & 0.37 & 2883 & 12.23 & 4.24 & 2883 & 21.42 & 7.43 \\
        GPBiCG\_dp  & 2883 & 2.09 & 0.72 & 2883 & 6.33 & 2.20 & 2883 & 12.67 & 4.39 \\
        GPBiCG\_d\_ILU(0)  & 1871 & 149.19 & 79.74 & 1162 & 100.73 & 86.68 & 896 & 84.30 & 94.09 \\ \hline
    \end{tabular}}
    \end{center}
\end{table}

Compared to the real-valued case, solving complex linear systems takes more time, consistent with our previous benchmark study \cite{kouya_iceet2024}. Furthermore, the cost per iteration when using ILU(0) preconditioning is significantly higher.

We also benchmarked arbitrary-precision arithmetic using MPFR and MPC with a 256-bit mantissa to observe the impact on convergence. Results indicate that like the multi-component case, no convergence is achieved without ILU(0). The effect of mixed-precision in these settings is summarized in Table~\ref{tab:mpfr}.

\begin{table}[htb]
    \caption{MPFR and MPC 256-bit benchmark results}
    \label{tab:mpfr}
    \begin{center}
    {\small \begin{tabular}{|c|rrr|rrr|} \hline
MPFR and MPC 256bits &  \multicolumn{3}{c|}{mcfe} &  \multicolumn{3}{c|}{dwg961b} \\ 
Algorithm & \multicolumn{1}{c|}{[1]} & \multicolumn{1}{c|}{[2]} & \multicolumn{1}{c|}{[2]/[1]} & \multicolumn{1}{c|}{[1]} & \multicolumn{1}{c|}{[2]} & \multicolumn{1}{c|}{[2]/[1]} \\ \hline
BiCG & 2295 & 11.93 & 5.20 & 2883 & 7.60 & 2.64 \\
BiCG\_p & 2295 & 7.31 & 3.19 & 2883 & 1.43 & 0.50 \\
BiCG\_ILU(0)  & 13 & 1.10 & 84.46 & 401 & 0.90 & 2.25 \\
BiCG\_d & 2295 & 17.70 & 7.71 & 2883 & 13.37 & 4.64 \\
BiCG\_dp  & 2295 & 7.33 & 3.20 & 2883 & 1.68 & 0.58 \\
BiCG\_d\_ILU(0)  & 13 & 0.86 & 66.31 & 2883 & 72.52 & 25.15 \\ \hline
CGS  & 2295 & 11.84 & 5.16 & 2883 & 9.78 & 3.39 \\
CGS\_p  & 2295 & 5.47 & 2.38 & 2883 & 1.34 & 0.46 \\
CGS\_ILU(0) & 9 & 0.82 & 91.44 & 468 & 0.59 & 1.26 \\
CGS\_d  & 2295 & 17.60 & 7.67 & 2883 & 17.80 & 6.17 \\
CGS\_dp & 2295 & 5.33 & 2.32 & 2883 & 1.68 & 0.58 \\
CGS\_d\_ILU(0) & 9 & 0.63 & 70.00 & 2087 & 0.47 & 0.22 \\ \hline
BiCGSTAB & 2295 & 12.12 & 5.28 & 2883 & 13.40 & 4.65 \\
BiCGSTAB\_p & 2295 & 5.74 & 2.50 & 2883 & 2.00 & 0.69 \\
BiCGSTAB\_ILU(0)  & 9 & 0.79 & 87.22 & 716 & 0.59 & 0.83 \\
BiCGSTAB\_d & 2295 & 17.88 & 7.79 & 2883 & 24.65 & 8.55 \\
BiCGSTAB\_dp  & 2295 & 5.53 & 2.41 & 2883 & 2.49 & 0.86 \\
BiCGSTAB\_d\_ILU(0)  & 9 & 0.60 & 66.89 & 2883 & 0.47 & 0.16 \\ \hline
GPBiCG & 2295 & 13.82 & 6.02 & 2883 & 17.60 & 6.10 \\
GPBiCG\_p & 2295 & 7.30 & 3.18 & 2883 & 3.40 & 1.18 \\
GPBiCG\_ILU(0)  & 8 & 0.76 & 95.00 & 674 & 0.76 & 1.13 \\
GPBiCG\_d & 2295 & 19.97 & 8.70 & 2883 & 31.06 & 10.77 \\
GPBiCG\_dp  & 2295 & 7.50 & 3.27 & 2883 & 3.96 & 1.37 \\
GPBiCG\_d\_ILU(0)  & 8 & 0.58 & 72.50 & 2883 & 0.49 & 0.17 \\ \hline
    \end{tabular}}
    \end{center}
\end{table}

Finally, Table~\ref{tab:ratio} summarizes the reduction ratios in computation time when using ILU(0) and mixed-precision ILU(0) substitutions.

\begin{table}[htb]
\begin{center}
    \caption{Computation time ratio: real vs. complex systems}\label{tab:ratio}
    {\small \begin{tabular}{|c|rrr|r|rrr|r|}\hline
    Ratio    & \multicolumn{4}{|c|}{mcfe} & \multicolumn{4}{c|}{dwg961b} \\ 
    Algorithm & \multicolumn{1}{c}{DD} & \multicolumn{1}{c}{TD} & \multicolumn{1}{c|}{QD} & \multicolumn{1}{c|}{MPFR} & \multicolumn{1}{c}{DD} & \multicolumn{1}{c}{TD} & \multicolumn{1}{c|}{QD} & \multicolumn{1}{c|}{MPC} \\ \hline
       ILU(0)/BiCG & 456.4 & 43.5 & 37.7 & 16.2 & 960.1 & 106.2 & 102.8 & 15.1 \\
       ILU(0)\_d/BiCG\_d & 476.7 & 42.4 & 37.2 & 8.6 & 313.6 & 29.5 & 42.8 & 8.2 \\
       ILU(0)/ILU(0)\_d & 1.0 & 1.1 & 1.1 & 1.3 & 3.2 & 3.7 & 2.5 & 1.2 \\ \hline
       ILU(0)/CGS & 538.0 & 45.0 & 37.8 & 17.7 & 863.5 & 94.3 & 86.0 & 14.5 \\
       ILU(0)\_d/CGS\_d & 515.4 & 43.5 & 35.7 & 9.1 & 284.3 & 25.8 & 17.7 & 7.8 \\
       ILU(0)/ILU(0)\_d & 1.1 & 1.1 & 1.1 & 1.3 & 3.2 & 3.8 & 5.0 & 1.2 \\ \hline
       ILU(0)/BiCGSTAB & 526.5 & 41.8 & 34.0 & 16.5 & 773.7 & 88.8 & 74.8 & 13.5 \\
       ILU(0)\_d/BiCGSTAB\_d & 505.7 & 40.4 & 32.2 & 8.6 & 258.0 & 24.3 & 15.6 & 7.5 \\
       ILU(0)/ILU(0)\_d & 1.1 & 1.1 & 1.0 & 1.3 & 3.1 & 3.8 & 4.9 & 1.2 \\ \hline
       ILU(0)/GPBiCG & 500.9 & 40.6 & 33.8 & 15.8 & 644.8 & 74.6 & 74.8 & 10.5 \\
       ILU(0)\_d/GPBiCG\_d & 481.5 & 40.4 & 32.2 & 8.3 & 212.9 & 20.4 & 12.7 & 6.4 \\
       ILU(0)/ILU(0)\_d & 1.1 & 1.1 & 1.1 & 1.3 & 3.1 & 3.7 & 4.8 & 1.2 \\ \hline
       \end{tabular}}
    \end{center}
\end{table}

The results indicate the following for real-valued systems:
\begin{itemize}
  \item ILU(0) preconditioning increases total time by factors of 456.4--538.0 (DD), 40.6--45.0 (TD), and 33.8--37.8 (QD).
  \item Mixed-precision substitutions reduce time by -4.5--5.3\% compared to standard ILU(0).
  \item For MPFR 256-bit, mixed-precision yields up to 47.1--48.5\% reduction in ILU(0) cost.
\end{itemize}

In complex systems, mixed-precision preconditioning significantly improves performance, especially for multi-component formats. However, in MPC-based arbitrary precision, the effect is less pronounced.
\begin{itemize}
    \item ILU(0) preconditioning increases total time by factors of 644.8--960.1 (DD), 74.6--106.2 (TD), and 74.8--102.8 (QD).
    \item Mixed-precision substitutions reduce time by 58.4--83.1\% compared to standard ILU(0).
    \item For MPFR 256-bit, mixed-precision yields up to 39.2--46.0\% reduction in ILU(0) cost.
  \end{itemize}

\section{Conclusion and Future Work}

The benchmark results presented above indicate that, for cases where preconditioning is essential, the computational cost of ILU(0) increases with higher arithmetic precision. Specifically, DD precision incurs a 456x to 960x overhead, TD incurs 41x to 106x, and QD incurs 34x to 103x compared to un-preconditioned runs. In such cases, the use of mixed-precision SpMV and mixed-precision forward/backward substitution provides a modest reduction in computational cost, typically around 58.4--83.1\% for multi-component types in complex Krylov subspace methods.

In the MPFR 256-bit precision tests, the benefit of using mixed-precision substitution was even more evident, reducing ILU(0) computation time by approximately 39.2--48.5\%. This shows that even a naive mixed-precision implementation can yield significant performance improvements, especially in high-precision environments.

Furthermore, for complex-valued systems, the benefit of mixed-precision is more pronounced when using multi-component types (DD, TD, QD). In contrast, the MPC-based arbitrary-precision results demonstrated minimal variation, likely due to the already high cost of complex arithmetic in those settings.

Currently, our ILU(0) implementation is not optimized or parallelized. The ILU factorization itself and the forward/backward substitutions are implemented sequentially, leading to significant computational overhead. To further improve performance, especially for practical use in high-precision environments, future work should include:
\begin{itemize}
  \item SIMD optimization and OpenMP parallelization of ILU(0) routines
  \item Improved memory access patterns and cache-aware implementations
  \item Enhanced support for mixed precision preconditioning techniques
  \item Benchmarking on wider hardware platforms including Arm NEON, WASM, and SSE
\end{itemize}

\section*{Acknowledgment}

This research was supported by JSPS KAKENHI Grant Number 23K11127.


\end{document}